\documentclass[french]{article} 
\usepackage[latin1]{inputenc}
\usepackage{fontenc}
\usepackage{amsmath}
\usepackage{amssymb}
\usepackage{graphicx}
\usepackage{epsfig}
\usepackage{babel}
\newtheorem{theorem}{Th\'eor\`eme} 
\newtheorem{definition}{D\'efinition}

\def\epsilon{\varepsilon}
\def\R{{\mathbb R}} 
\def\C{{\mathbb C}} 
\def\Z{{\mathbb Z}} 
 
\def\N{{\mathbb N}}
\def\v{\text{vol}}

\def\C{{\tilde{C}}} 

\def\epsilon{{\varepsilon}} 
 
\def\h{{\widehat{h}}} 
\def\v{{\bf v}} 
\def\B{{\cal B}_{st}} 
\def\T{{\mathbb T}} 
\def\J{{\cal J}}

\title{Sur la forme de la boule unité de la norme stable unidimensionnelle}
\author{Ivan BABENKO \footnote{{UMR 5149, Institut de Math\'ematiques et de Modélisation
de  Montpellier, Universit\'e Montpellier II Case Courrier 051 - Place Eug\`ene Bataillon 34095
  Montpellier CEDEX 5, France e-mail: babenko@math.univ-montp2.fr}},
 Florent BALACHEFF\footnote{UMR 5149, Institut de Math\'ematiques et de Modélisation
de  Montpellier, Universit\'e Montpellier II  Case Courrier 051 - Place Eug\`ene Bataillon 34095
  Montpellier CEDEX 5, France e-mail: balachef@math.univ-montp2.fr}} 

\begin{document}

\maketitle

\begin{abstract}

Pour un polyèdre riemannien, nous étudions les formes apparaissant comme boule unité de la norme stable
unidimensionnelle (boule stable). Dans le cas d'un polyèdre riemannien unidimensionnel (graphe), nous montrons que la boule stable est un polytope dont les sommets
sont complètement décrits par la combinatoire du graphe. Nous étudions ensuite les formes réalisables comme boule stable de variétés riemanniennes de dimension plus grande que trois. Nous montrons
que, pour une variété riemannienne $(M,g)$ fixée, une large classe de polytopes
peut apparaître comme boule stable d'une métrique dans la classe conforme de
$g$. Nous utilisons pour cela une technique polyèdrale.

\end{abstract}


\noindent {\it Mots clefs :} variété riemannienne, poly\`edre riemannien,
norme stable, boule unité. 
\bigskip

\section{Introduction}

\subsection{Enoncé des résultats}

Pour un poly\`edre fini riemannien $(M,g)$ de classe $C^0$, 
l'homologie réelle possède une norme naturelle appelée {\it masse} 
ou {\it norme stable} (voir \cite{fede}, \cite {GLP}). Cette norme 
est particulièrement intéressante pour l'homologie de dimension $1$ :
 elle contrôle le comportement à l'infini de la géométrie relevée sur 
le revêtement homologique correspondant. Nous rappelons ici la plus 
géométrique des définitions en renvoyant le lecteur à \cite{fede} 
pour les différentes approches et 
la démonstration des équivalences entre celles-ci.

\begin{definition}
Soit $v \in H_1(M,\R)$ une classe entière. On pose
$$
\|v\|_{g}=\lim_{n \rightarrow + \infty} \frac{L_g(\gamma_n)}{n}, 
$$
où $\gamma_n$ est la plus petite courbe fermée réalisant 
la classe $n.v$ et $L_g$ désigne la longueur respectivement à $g$.
Cela définit par continuité une norme sur $H_1(M,\R)$ 
appelée {\it norme stable}. On pose alors :
$$
\B(g)=\{u \in H_1(M,\R) \mid \|u\|_g\leq 1\}.
$$

\end{definition}

Un point $s$ de $H_1(M,\R)$ est dit {\it de direction rationnel} 
si il existe une classe entière $v$ telle que 
$s \in <v>$ la droite engendrée par $v$. 
Dans le cas contraire, le point est dit {\it de direction irrationnelle}.

Les vari\'et\'es riemanniennes et les polyèdres unidimensionnels (ou graphes) constituent deux sous-classes naturelles de poly\`edres riemanniens pour lesquelles nous allons \'etudier les formes réalisables comme boule unit\'e de la norme stable associée à une métrique.

Nous commen\c cons par étudier la classe des graphes finis. La m\'etrique naturelle d'un graphe pond\'er\'e est celle induite par la métrique
sur ses ar\^etes vues comme intervalles de longueur donnée par le poids.
Comme nous le démontrons dans la section 2.1, la boule unit\'e
de la norme stable d'un graphe pondéré est toujours un polytope. On peut d\'ecrire complètement 
la forme de ce polytope en terme de la combinatoire du graphe.
Pour cela, étant donné un graphe $G$, un chemin ferm\'e orienté parcourant tous ses sommets
une seule fois est dit {\it circuit simple orienté}. Deux circuits simples qui co\"incident géométriquement et sont de même orientation, mais qui diffèrent par leur point initial sont identifiés. On consid\`ere l'espace vectoriel 
${\cal{C}}(G,\R)$ engendr\'e par les arêtes de $G$ après le choix d'une orientation arbitraire. 
A chaque circuit simple orienté est associé un vecteur de ${\cal{C}}(G,\R)$ - la somme (avec les signes n\'ecessaires) des ar\^etes formant ce circuit.
Les vecteurs ainsi construits sont de nouveau appelés circuits simples orientés. Notons qu'à un circuit géométrique simple correspond deux circuits simples orientés.

\vskip7pt
\noindent
{\bf Th\'eor\`eme {\bf A}}  {\it Soit $(G,w)$ un graphe pondéré. On note
$\{C_j\}_{j \in J}$ l'ensemble de ses circuits simples orient\'es.

Alors la boule unit\'e de la norme stable pour la métrique $w$ dans $H_1(G, \R)$
co\"{\i}ncide avec l'enveloppe convexe dans ${\cal{C}}(G,\R)$ 
des vecteurs $\{C_j/\|C_j\|_w \}_{j \in J}$.}
 
\vskip7pt

A premier nombre de Betti $b$ fixé, le nombre de circuits simples orientés d'un graphe admet une majoration évidente : 

\vskip7pt
\noindent
{\bf Corollaire {\bf A}}  {\it Soit $(G,w)$ un graphe pondéré 
de premier nombre de Betti $b$.  
Alors la boule unit\'e de sa norme stable est un
polytope $b$-dimensionnel dont le nombre de sommets est majoré par $2(2^b - 1)$.}

\vskip7pt

Etant fixée une variété fermée $M$, nous allons nous intéresser 
maintenant aux  normes  qui peuvent être réalisées comme norme stable 
d'une métrique lisse ({\it i.e} de classe $C^\infty$).

\vskip7pt
\noindent
{\bf Th\'eor\`eme {\bf B}}  {\it Soit $(M, g)$ une variété fermée riemannienne lisse  de 
dimension $m \geq 3$ et de premier nombre de Betti ${b \geq 1}$ .
On consid\`ere  un polytope fini convexe $K$ de 
$H_1(M, {\R})$, à symétrie centrale, tel que son int\'erieur soit non vide
et les directions de ses sommets  soient rationnelles.
Alors il existe une m\'etrique $g'$ lisse sur $M$ conforme \`a $g$
telle que}
$$ \B(g') = K. $$
\vskip7pt

On en déduit immédiatement le corollaire suivant :

\vskip7pt
\noindent
{\bf Corollaire {\bf B}}  {\it Soient $(M, g)$ une vari\'et\'e riemannienne de 
dimension $m \geq 3$ et de premier nombre de Betti $b\geq 1$ et
$\| \cdot \|$ une norme sur l'espace vectoriel $H_1(M, \R)$. 
Pour tout $\varepsilon > 0$, il existe une m\'etrique 
$g(\varepsilon)$ conforme \`a $g$ telle que}
$$ 
\| \cdot \| \leq \| \cdot \|_{g(\varepsilon)}
\leq  (1 + \varepsilon)\| \cdot \| .
$$

\subsection{Historique}

Nous récapitulons ici l'ensemble des résultats traitant de la question suivante : 
étant fixée une variété $M$, quelles normes peuvent être réalisées 
comme norme stable associée à une métrique $g$ ?

\bigskip
\noindent {\bf Cas du tore bidimensionnel}. 
M.Morse \cite{mors} a tout d'abord remarqué  que la norme stable 
du tore bidimensionnel $\T^2$ associée à une métrique lisse est 
strictement convexe. V. Bangert \cite{bang2} a ensuite montré le résultat suivant :
\begin{theorem}
Soit $(\T^2,g)$ un tore riemannien lisse. La norme stable est \linebreak 
différentiable dans les points de direction irrationnelle et 
est différentiable dans les points de direction rationnelle 
si et seulement si le tore est feuilletté par les géodésiques 
minimisantes représentant l'élément correspondant de $H_1(\T^2,\Z)$.
\end{theorem}

\bigskip
\noindent {\bf Cas des surfaces}.  On note $\Sigma_h$ la surface fermée
orientable de genre $h$. On se fixe une métrique $g$ lisse. 
Un {\it sous-espace affine support} de $\B(g)$ est 
un sous-espace affine $H$ tel que $H \cap \B(g) = H \cap \partial \B(g)$. 
On dit d'un sous-ensemble $F$  de $\partial \B(g)$ qu'elle est une {\it face plate} 
si il existe un sous-espace affine support $H$ tel que $F=H \cap \partial \B(g)$. 
La dimension de $F$ est alors définie comme la dimension de l'espace affine $A$ qu'elle engendre et l'intérieur de $F$ est son intérieur dans $A$.
Enfin une face plate est {\it rationnelle} si elle co\"incide avec 
l'enveloppe convexe de points de $\partial \B(g)$ de direction rationnelle. 
Dans \cite{mass}, D.Massart montre que l'ensemble des faces contenant 
un point de $\partial \B(g)$ dans leur intérieur peut être ordonné 
par inclusion et admet un maximum unique. 
Il montre également le résultat suivant :

\begin{theorem}
Soit $g$ une métrique lisse sur $\Sigma_h$ avec $h\geq 2$. 

\noindent 1. La dimension d'une face plate est majorée par $h-1$. 

\noindent 2. Tout point de direction rationnelle est contenu dans 
une face plate rationnelle de dimension $h-1$.

\noindent 3. En un point $v$ de $\partial \B(g)$ de direction rationnelle, 
la norme stable est différentiable seulement dans les directions tangentes 
à $F$, où $F$ est la face plate rationnelle maximale 
contenant $v$ dans son intérieur.
\end{theorem} 

On en déduit immédiatement que la norme stable d'une surface 
$(\Sigma_h,g)$ de genre $h \geq 2$ n'est donc ni lisse, ni strictement 
convexe, ni définie par un polytope. Remarquons également qu'il existe 
une infinité de points exposés de direction rationnelle pour de telles 
surfaces (un point est dit {\it exposé} si la face rationnelle maximale 
le contenant dans son intérieur est réduite à ce point) (voir \cite{mass2}).

\medskip

Dans \cite{mass3}, D.Massart relie la dérivabilité de la norme stable 
en une classe d´homologie au degré d'irrationalité de cette classe. 
Plus précisément, une classe d´homologie $v$ est dite {\it $k$-irrationnelle} 
si $k$ est la dimension du plus petit sous-espace de $H_1(\Sigma_h,\R)$ 
engendré par des classes entières et contenant $v$. 

\begin{theorem}
Soit $g$ une métrique lisse sur $\Sigma_h$ avec $h\geq 1$. 
En une classe $v$ d'homologie $k$-irrationnelle, 
la norme stable est différentiable dans au moins 
$k-1$ directions non radiales.
\end{theorem}

\bigskip
\noindent {\bf Cas des dimensions supérieures}. 
G.A.Hedlund a construit pour le tore de dimension 3 une métrique Riemannienne \cite{hedl} pour laquelle la norme stable est donnée par un octaèdre régulier. Cet exemple est fondamental, puisqu'il montre que la situation en dimension plus grande que 3 est radicalement différente de celle en dimension 2 : les polytopes peuvent apparaître comme boule unité d'une norme stable associée à une métrique.

\medskip

D'autre part, D.Burago, S.Ivanov et B.Kleiner \cite{buivkl} ont démontré le résultat suivant, qui généralise le résultat de V.Bangert :

\begin{theorem}
Soit $(\T^n,g)$ un tore riemannien de classe $C^3$ de dimension 
$n\geq 2$. Soit $v$ un point de direction irrationnelle. 
Alors la norme stable est différentiable en $v$ dans 
au moins une direction non radiale.
\end{theorem} 

Dans ce même papier, les auteurs montrent que, 
pour tout entier $k$, il existe un entier $n$ 
et une métrique de classe $C^k$ sur $\T^n$ 
telle que pour presque tous les points de 
direction irrationnelle $v$, la norme stable 
soit non différentiable en $v$. Cela montre 
qu'une généralisation complète du théorème de 
V.Bangert concernant les directions irrationnelles est impossible.

\bigskip

Les auteurs expriment leurs remerciements à D.Massart pour d'instructives et de stimulantes conversations.

\bigskip

Le papier est organisé de la manière suivante : les sections 2 et 3 correspondent respectivement aux démonstrations des théorèmes {\bf A} et {\bf B}.

\section{Norme stable des graphes}

\subsection{Rappels sur les graphes}

Commen\c cons par quelques définitions. 
Un {\it graphe pond\'er\'e} est une paire $(G,w)$ o\`u $G=(V,E)$ 
est un multigraphe fini connexe non orienté et $w$ 
est une {\it fonction poids} sur les ar\^etes 
$w:E \rightarrow \mathbb R_+$. 
Un multigraphe désigne un graphe dans lequel 
on autorise les arêtes multiples et les boucles. 
On appelle $w(e)$ le {\it poids} d'une ar\^ete. 
Tout graphe est naturellement identifi\'e \`a 
un graphe pond\'er\'e dans lequel le poids 
de chaque ar\^ete vaut 1. 

\medskip

On se fixe pour la suite un graphe pondéré $(G,w)$ 
de premier nombre de Betti $b$. On choisit pour chaque arête une orientation arbitraire, notons ces arêtes $\{e_i\}_{i=1}^k$. L'espace vectoriel ${\cal{C}}(G,\R)$ engendré par les arêtes orientées $\{e_i\}_{i=1}^k$ co\"incide avec l'espace des chaînes simpliciales du complexe simplicial $G$ :
$$
{\cal{C}}(G,\R) =\{ \sum_{i=1}^k a_i.e_i \mid  a_i \in \R \text{ pour } i=1, \ldots, k\}.
$$
Comme les graphes n'ont pas de cellule pour 
les dimensions plus grandes que $2$, 
l'homologie de $G$ de dimension $1$ à coefficients réels 
$H_1(G,\R)$ est plongée naturellement dans 
${\cal{C}}(G)$ 
comme un sous-espace vectoriel de dimension $b$. 
L'homologie de $G$ de dimension $1$ à coefficients 
entiers $H_1(G,\Z)$, en l'absence de torsion dans ce 
cadre unidimensionnel, constitue un réseau du sous-espace 
$H_1(G,\R)$ (comparer avec \cite{bacharnag}).

\medskip

 Pour $u=\sum_{i=1}^k u_i.e_i \in {\cal{C}}(G,\R)$, on note 
$$| u |_{w,1}=\sum_{i=1}^k w_i |u_i|, \eqno (2.1)
$$
où $w_i=w(e_i)$ pour $i=1,\ldots,k$. On voit facilement que cette norme co\"incide avec 
la norme stable $\| \cdot \|_w$ sur $H_1(G,\R)$. En effet, la norme stable de $v \in H_1(G,\R)$ est donnée par la formule
$$
\|v\|_{w} = \inf \{\sum_{i=1}^s |\alpha_i| w(\sigma_i) \mid v = \sum _{i=1}^s \alpha_i [\sigma_i] , \alpha_i \in \R \text{ et } \sigma_i \in E \}.
$$
Cette définition est équivalente à la définition 1 (voir \cite{fede}).

La boule stable $\B(G, w)$ 
est donc l'intersection 
de la boule unité de la norme $| \cdot |_{w,1}$ 
dans l'espace des arêtes ${\cal{C}} (G,\R)$ avec le sous-espace vectoriel 
de dimension $b$ image par le plongement naturel 
de l'homologie réelle de $G$ de dimension $1$.

\subsection{Démonstration du théorème A}

La d\'emonstration du th\'eor\`eme {\bf A} est une cons\'equence
imm\'ediate des lemmes suivants.

\vskip7pt
\noindent
{\bf Lemme A1} {\it Tout circuit simple orienté de $G$, identifié au
vecteur correspondant de ${\cal{C}}(G,\R)$, est proportionnel à un sommet de $\B(G, w)$. Le facteur de proportionnalité est exactement la longueur de ce circuit.}
\vskip7pt

\noindent 
{\bf D\'emonstration.} Notons $B_{w,1}$ la boule unit\'e pour la
norme (2.1) dans ${\cal{C}}(G,\R)$. Comme
$$
 \B(G, w) = B_{w,1} \bigcap H_1(G,\R),
 $$
 les sommets de $\B(G, w)$   
sont donnés par les points d'intersection de $H_1(G,\R)$ avec l'intérieur des
faces de $B_{w,1}$ de codimension plus grande ou \'egale \`a $b$ dans le cas où cette intersection est réduite à un point.
Pour $C$ un circuit simple orienté, on note 
$$
C = \mathop{\sum} \limits_{i \in I(C)}
\varepsilon_i e_i  \eqno (2.2) 
$$
son d\'eveloppement dans la base des ar\^etes $\{e_i\}_{i=1}^k$ (ici, $\varepsilon_i = \pm 1$). Il est \'evident que le nombre 
$|I(C)|$ des ar\^etes dans (2.2) n'exc\`ede pas $k - b + 1$.
On consid\`ere la face $|I(C)|-1$ dimensionnelle $F(C)$ de $B_{w,1}$ contenant les vecteurs
$$ 
{\varepsilon_i \over w(e_i)}e_i , \hskip5pt i \in I(C) .
$$
On voit facilement que 
$$ 
X= {1\over \|C\|_w} C = {1 \over \mathop{\sum}\limits_{i \in I(C)} w(e_i)}
\mathop{\sum}\limits_{i \in I(C)} \varepsilon_i e_i   
$$
est un point de $\mbox{int} (conv(\{{\varepsilon_i \over w(e_i)}e_i, i \in I(C)\}))
$ où $conv (A)$ désigne l'enveloppe convexe d'un ensemble fini de points $A$ et $\mbox{int} (B)$ désigne l'intérieur d'un ensemble $B$.

On montre alors que l'intersection 
$H_1(G,\R) \bigcap \mbox{int}(F(C))$ est réduite à ce point. En effet, si $v$ est un autre point de $H_1(G,\R) \bigcap \mbox{int}(F(C))$,  les points $v$ et $X$ définissent un segment contenu dans l'intersection. Ceci implique qu'il existe des points d'intersection entre $H_1(G,\R)$ et $F(C)$ dans un voisinage arbitraire de $X$. Autrement dit, il existe des points d'intersection (différents de $X$) appartenant à l'enveloppe convexe des $\{{\varepsilon_i \over w(e_i)}e_i, i \in I(C)\}$. Si $u$ est un tel point, on  peut écrire $u$ comme une combinaison linéaire des $\{{\varepsilon_i \over w(e_i)}e_i, i \in I(C)\}$ dont les coefficents sont tous non nuls :
$$
u=\sum_{i \in I(C)} \alpha_i {\varepsilon_i \over w(e_i)}e_i.
$$
 On complète alors le vecteur $C$ en une base de $H_1(G,\R)$ formée de circuits simples $\{C_1=C,C_2,\ldots,C_b\}$ de sorte que, pour chaque $i$ dans $\{1,\ldots,b\}$, il existe une arête $f_i$ contenue dans $C_i$ et n'appartenant pas aux $C_j$ pour $j\neq i$. L'analyse de la décomposition du vecteur $u$ dans cette base nous montre que le vecteur $u$ est nécessairement proportionnel à $C$. D'où une contradiction.

Maintenant, comme $\mbox{dim} (F(C)) \leq k - b $, le vecteur $C/ \|C\|_w$ est bien un sommet de $\B(G, w)$. Ceci ach\`eve la d\'emonstration.   \hspace{\stretch{1}} $\Box$\\

\vskip7pt
\noindent
{\bf Lemme A2} {\it Soient $C_1$ et $C_2$ deux circuits simples orientés. 
Alors il existe des circuits simples orientés
$\{D_j\}_{j \in J}$ (non uniquement d\'efinis) pour lesquels chaque
ar\^ete n'est parcourue que dans un seul sens et tels que dans  $H_1(G,\R)$ :}
$$[ C_1 + C_2] = \mathop{\sum}\limits_{j \in J} [D_j].$$
\vskip7pt

\noindent 
{\bf D\'emonstration.} On note $\{e_i\}_{i \in I}$  les ar\^etes de $G$ figurant
dans $C_1$ et $C_2$ avec des directions oppos\'ees. On \^ote les ar\^etes 
$\{e_i\}_{i \in I}$ de la reunion $C_1 \bigcup C_2$.
Notons $C$ la courbe obtenue. Elle poss\`ede une orientation
induite par celles de  $C_1$ et de $C_2$. En suivant cette orientation,
$C$ se sépare naturellement en quelques courbes ferm\'ees et
orient\'ees $\{P_l\}_{l \in L}$. Cette repr\'esentation
préserve la classe d'homologie :  
  $$
[C_1 + C_2] =  \mathop{\sum}\limits_{l \in L} [P_l] .   
$$
Chaque courbe $P_l$ n'est pas, en g\'en\'eral, un circuit
simple car elle peut avoir des auto-intersections : des sommets,
ou bien des ar\^etes parcourues plusieurs fois dans la m\^eme 
direction. On part d'un point d'auto-intersection et on parcourt $P_l$ en suivant l'orientation. En revenant au point de
départ pour la premi\`ere fois, nous coupons $P_l$ en deux
courbes ferm\'ees et orient\'ees. Chacune de ces deux nouvelles 
courbes a moins d'auto-intersections que la courbe initiale
$P_l$. En r\'ep\'etant suffisamment ce proc\'ed\'e, nous
obtenons un certain nombre de circuits simples qui engendrent
la classe $[P_l]$. On applique ce proc\'edé pour chaque $P_l , \hskip5pt l \in L$. L'ensemble des circuits simples
obtenus est noté $\{D_j\}_{j \in J}$. Pour achever la d\'emonstration, il reste \`a remarquer que, par construction, les circuits simples orientés construits n'ont pas d'ar\^ete commune parcourue dans des directions oppos\'ees.  \hspace{\stretch{1}} $\Box$\\

\vskip7pt
\noindent
{\bf Lemme A3} {\it Pour chaque classe enti\`ere $a \in H_1(G,\Z)$,
il existe des circuits simples $\{C_s\}_{s \in S}$ (non uniquement
d\'efinis) tels que }
$$ a = \mathop{\sum}\limits_{s \in S}[C_s] \hskip5pt \mbox{ et }
\hskip5pt \|a\|_w = \mathop{\sum}\limits_{s \in S} \|C_s\|_w .$$
\vskip7pt

\noindent 
{\bf D\'emonstration.} Comme les circuits simples orientés engendrent
 $ H_1(G,\Z)$, nous pouvons présenter une classe fixée  $a \in H_1(G,\Z)$ comme une somme de circuits simples orientés :
$$ 
a = \mathop{\sum}\limits_{r \in R}[D_r] . \eqno(2.3) 
$$
Remarquons que la repr\'esentation (2.3) n'est pas uniquement d\'efinie
et que certains circuits figurent, en g\'en\'eral, plusieurs fois 
dans cette somme. En appliquant syst\'ematiquement le lemme {\bf A2} sur
chaque paire de circuits $(D_{r_1},D_{r_2})$ ayant des ar\^etes
en commun parcourues dans des directions oppos\'ees, nous arrivons
au bout de ce proc\'ed\'e itératif sur une  nouvelle représentation par des circuits
$\{C_s\}_{s \in S}$ v\'erifiant le lemme. \hspace{\stretch{1}} $\Box$\\

\vskip7pt
\noindent
{\bf Lemme A4} {\it Soit  $(G,w)$ un graphe pond\'er\'e dont la
fonction poids $w$ est à valeurs rationnelles.
Alors chaque sommet de  $\B(G, w)$ est proportionnel \`a
un circuit simple orienté de $G$. }
\vskip7pt
\noindent 
{\bf D\'emonstration.} On consid\`ere un sommet $X$ de $\B(G, w)$.
On sait que c'est l'unique point d'intersection de $H_1(G,\R)$ et
d'une face de  $B_{w,1}$ de codimension plus grande ou \'egale \`a
$b$. Comme les poids $\{w(e_i) ,\hskip5pt 1 \leq i \leq k\}$ sont rationnels,
les coordonn\'ees des sommets de $B_{w,1}$ sont rationnelles.
Le sous-espace $H_1(G,\R)$ est engendr\'e par des vecteurs dont les coordonn\'ees
sont enti\`eres. Ceci implique que les coordonn\'ees de $X$ sont rationnelles.
On constate donc que $X = \lambda a$, o\`u $a \in H_1(G,\Z)$ est un vecteur
indivisible et $\lambda$ est un facteur rationnel positif. On a \'egalement
$$ 1 = \|X\|_w = \lambda \|a\|_w  . \eqno (2.4) $$   
D\'ecomposons $a$ en circuits simples orientés selon le lemme {\bf A3} :
$$ a = \mathop{\sum}\limits_{s \in S}[C_s]. $$
On tire alors du lemme {\bf A3} et de (2.4) l'\'egalit\'e suivante :
$$ X = {1 \over \mathop{\sum}\limits_{s \in S} \|C_s\|_w }
\mathop{\sum}\limits_{s \in S} \Biggl{(} \|C_s\|_w \Biggl{(} {1 \over \|C_s\|_w } 
[C_s].   \Biggr{)} \Biggr{)}  \eqno (2.5)$$
D'apr\`es le lemme {\bf A1}, les vecteurs $[C_s]/\|C_s\|_w$ pour $s \in S$ sont des sommets de $\B(G, w)$ et (2.5) implique \`a son tour que
$X$ est un point de l'intérieur de l'enveloppe convexe de ces points. Comme $X$
est un sommet, on obtient donc que $S$ contient un seul circuit (la r\'ep\'etition
d'un unique circuit n'étant pas possible, puisque $a$ a \'et\'e choisi indivisible).
La d\'emonstration est ach\'ev\'ee. \hspace{\stretch{1}} $\Box$\\

Les lemmes {\bf A1} et {\bf A4} impliquent le th\'eor\`eme {\bf A} pour 
tout graphe pond\'er\'e dont la fonction poids est \`a valeurs rationnelles.
D'autre part, ces lemmes montrent \'egalement que les directions des sommets
de $\B(G, w)$ ($w$ étant \`a valeurs rationnelles)
sont uniquement d\'efinies et ne d\'ependent que des circuit simples. 
Ceci, par continuit\'e, implique le r\'esultat pour des poids quelconques.

\section{Polytope et classe conforme}

La d\'emonstration du th\'eor\`eme va découler de plusieurs
propositions techniques. On travaille dans la classe des métriques simpliciales. Ce sont des métriques lisses sur chaque simplexe d'une triangulation lisse avec une condition naturelle de co\"incidence sur les faces communes. Pour plus de détails, le lecteur pourra consulter \cite{babe}.

\subsection{Préparatifs}

\vskip7pt
\noindent
{\bf Lemme B1} {\it  Soient $(X_i, h_i) , \hskip5pt i= 1,2$
deux variétés munies de métriques simpliciales et $f: X_1 \longrightarrow X_2$
une application simpliciale contractant les distances. On note $f_*$ l'application induite sur les homologies r\'eelles unidimensionnelles.
Alors
$$ f_*(\B(h_1)) \subset \B(h_2) . $$}
\vskip7pt

\noindent {\bf D\'emonstration.} Soit ${\bf a}$ un élément de $H_1(X_1, \Z)$ et supposons que $\gamma_n$ soit une courbe fermée lisse par morceaux réalisant la classe 
$n{\bf  a}$. Comme f contracte les distances, on a
$ l_{h_1}(\gamma_n) \geq l_{h_2}(f(\gamma_{n})) $, ce qui implique que 
$$\|{\bf a}\|_{h_1} \geq  
 \|f_*({\bf a})\|_{h_2} . \eqno(3.1)  $$
Par continuité de la norme, l'inégalité (3.1) est valable pour toute classe réelle. On en déduit le résultat. \hspace{\stretch{1}} $\Box$\\

\bigskip

On note $\{s_1,\ldots,s_k\}$ un sous-ensemble de sommets de $K$ tels que $$
K=conv_s\{s_1,\ldots,s_k\}
$$
 où $conv_s$ désigne l'enveloppe convexe du symétrisé d'un ensemble. Les sommets de $K$ étant rationnels, on peut trouver  pour tout $i=1\ldots k$ un nombre $l_i>0$ tel que
$$
s_i=l_i.\v_i,
$$
o\`u $\v_i$ est un vecteur entier indivisible.
On se fixe $k$ courbes de $M$ lisses, disjointes, fermées et simples notées $\gamma_1,\ldots,\gamma_k$ telles que $[\gamma_i]=\v_i$ pour $i=1,\ldots,k$. On les suppose paramétrisées par longueur d'arc respectivement à la métrique donnée $g$.

\vskip7pt
\noindent
{\bf Lemme B2} {\it Dans $(M,g)$, il existe $k$ voisinages tubulaires disjoints $\{U_i\}_{i=1}^k$
des courbes  $\{\gamma_i\}_{i=1}^k$ et une m\'etrique $g_1$ lisse conforme 
\`a $g$ tels que 

\noindent
1. $l_{g_1}(\gamma_i) = l_i^{-1} , \hskip5pt 1 \leq i \leq k$ ;

\noindent
2. les $U_i$ sont fibr\'es par des disques $g_1$-orthogonaux \`a $\gamma_i$
et les projections \linebreak $p_i: U_i \longrightarrow \gamma_i$
le long de ces disques contractent les $g_1$-distances.}
\vskip7pt

\noindent
{\bf D\'emonstration.} On se fixe $\varepsilon > 0$ suffisamment petit pour que les $\varepsilon$-voisinages tubulaires des
courbes $\gamma_i$ respectivement à la métrique $g$ soient disjoints et qu'ils soient fibr\'es par des disques
$g$-g\'eodesiques orthogonaux aux $\gamma_i$. Notons $U_i$ ces voisinages tubulaires. On choisit une fonction $\lambda$ 
strictement positive sur $M$ qui soit constante égale à $l_i . l_g(\gamma_i)^{-1} $
sur $U_i$. La metrique $\widehat{g} = {\lambda}^2g$
v\'erifie la propriété 1. et les disques fibrant $U_i$ sont toujours 
g\'eodesiques et ortogonaux aux $\gamma_i$ pour la métrique $\widehat{g}$.
Puisque  $\widehat{g}$ est conforme \`a $g$, et pour simplifier nos
futures notations, nous pouvons supposer que la m\'etrique de départ
v\'erifie la propri\'et\'e 1.

Soit $p_i: U_i \longrightarrow \gamma_i$ le fibr\'e normale de $\gamma_i$.
On consid\`ere un syst\`eme de coordonn\'ees semi-g\'eodesique
sur $U_i$ :
$$(s,x_2,x_3, ... ,x_m), \eqno(3.2) $$
ou $s=x_1$ est le param\'etre 
naturel le long de $\gamma_i$ et les $(x_2,x_3, ... ,x_m)=\overline{x}$ sont les coordonn\'ees
géodésiques sur les disques ortogonaux. Ce syst\`eme est bien d\'efini
sur chaque ouvert du type $p_i^{-1}([s-\delta, s+\delta] $ pour $\delta > 0$ suffisamment petit. Pour $q=(s,\overline{x}) \in U_i$, notons également ${\bf r}(q) = \mbox{dist}(q, \gamma_i)$. Alors ${\bf r}^2$ 
est une fonction lisse définie sur $U_i$.
Soient $\{g_{kl}(s,\overline{x})\}^m_{k,l = 1}$ les 
coefficients de $g$ dans (3.2). Nous avons par construction 
$$  g_{11}(s,\overline{0}) = 1 ;
\hskip5pt g_{1l}(s,\overline{0}) = 0 , \hskip5pt 2 \leq l \leq m . \eqno(3.3) $$
Tout vecteur tangent ${\bf v} \in T_q U_i$ s'écrit dans la base $\{{\partial \over \partial x_k}\}_{k=1}^m $ 
$$
{\bf v} = (\alpha, \overline{v}) ; \hskip5pt \mbox{o\`u} \hskip5pt
  \overline{v} = (v_2,v_3,...,v_m). 
$$ 
Nous pouvons calculer sa norme :
$$\|{\bf v} \|_g^2 = g_{11}(s,\overline{x}) \alpha^2 + 2\alpha \mathop{\sum} \limits_{k=2}^m
g_{1k}(s,\overline{x})v_k +  \| \overline{v} \|_{g'}^2 , \eqno(3.4) $$
o\`u $g'$ est la restriction de $g$ à l'espace tangent au disque
normal correpondant, sa matrice co\"incidant au bloc de $g$ dual \`a 
$g_{11}$.

Cherchons maintenant le facteur conforme d\'efini sur $U_i$ sous la
forme suivante :
$$ \lambda_i^2(q) = {\frac {1+a{\bf r}(q)^2} {g_{11}(q)} } ,\eqno(3.5) $$
o\`u $a$ est une constante que nous allons choisir afin de v\'erifier
la propriété 2. du lemme.  Sur chaque $U_i$, pour la m\'etrique $g_1 =  \lambda_i^2g$, on tire de (3.4) que :
$$\|{\bf v} \|_{g_1}^2 - \alpha^2 = \lambda_i^2 \Biggl( {\frac{ a
    g_{11}(s,\overline{x}) {\bf r}^2(s,\overline{x}) }
{1+a{\bf r}^2(s,\overline{x})  }}\alpha^2  + 
2\alpha \mathop{\sum} \limits_{k=2}^m
g_{1k}(s,\overline{x})v_k +  \| \overline{v}
\|_{g'}^2\Biggr). \eqno(3.6) $$
Compte-tenu de (3.3), on a 
$$| \mathop{\sum} \limits_{k=2}^m
g_{1k}v_k | \leq A{\bf r}(s,\overline{x}) \| \overline{v}
\|_{g'} $$
pour une constante $A > 0$ appropri\'ee.
Nous obtenons alors de (6) l'inégalit\'e suivante :
$$ \|{\bf v} \|_{g_1}^2 - \alpha^2 \geq 
 \lambda_i^2 \Biggl( {\frac{ a
    g_{11}{\bf r}^2 }
{1+a{\bf r}^2  }}\alpha^2 - 2A|\alpha|{\bf r} \| \overline{v}\|_{g'} +
 \| \overline{v} \|_{g'}^2\Biggr). \eqno(3.7)
$$

Nous avons dans les parenth\`eses une  forme quadratique
en les arguments $|\alpha|{\bf r} $ et $ \| \overline{v}\|_{g'}$,
qui est positive si
$${\frac{ ag_{11}}  {1+a{\bf r}^2}} >
A^2 .$$
Prenons $a = 2A^2$, alors 
 cette in\'egalit\'e est v\'erifi\'ee pour
${\bf r} = 0$  en accord avec (3.3). 
On en d\'eduit par continuit\'e qu'elle est v\'erifi\'ee sur
l'ouvert $U_i$ tout entier si l'on a choisi $\varepsilon$ suffisamment petit au départ. Finalement, pour le facteur conforme (3.5) 
avec $a = 2A^2$, on obtient de (3.7) 
$$\|{\bf v} \|_{g_1} \geq |\alpha| $$
et la propriété 2. est donc vérifiée.

Pour achever la d\'emonstration, il reste \`a prolonger
les facteurs locaux conformes $\{\lambda_i \}^n_{i=1} $
en une fonction lisse $\lambda$ strictement positive.  \hspace{\stretch{1}} $\Box$\\

\subsection{Construction d'une métrique de référence}

Nous allons commencer par construire une métrique de référence sur $M$ pour laquelle la boule unité de la norme stable est bien le polytope demandé. Cette métrique n'appartient en général pas à la classe conforme de la métrique de départ.

\vskip7pt
\noindent
{\bf Proposition B3} {\it 
Il existe une m\'etrique lisse $h$ sur $M$ tel que 

\noindent 
1. $\|{d\gamma_i(s) \over ds}\|_h = \|{d\gamma_i(s) \over ds}\|_{g_1}$ , 
o\`u $s$ est le param\`etre $g_1$-naturel le long des $\gamma_i$ ;

\noindent 
2. ${\cal B}_{st}(M,h)=K$.}
\vskip7pt

\noindent {\bf D\'emonstration.}
La démonstration se déroule en deux étapes. Dans la première, on montre que l'on peut se ramener à construire une métrique simpliciale sur $M$ réalisant le polytope $K$ comme boule stable, et dans une seconde étape, on explicite la construction de cette métrique simpliciale.

\medskip

\noindent
{\bf Lemme B4} {\it 
Supposons qu'il existe une métrique $\h$ sur $M$ simpliciale, telle que :

\noindent 
1. $\B(\h)=K.$

\noindent 
2. $\exists \epsilon >0$ tels que la métrique $\h$ soit lisse sur les $\epsilon$-voisinages tubulaires $U_\epsilon(\gamma_i)$ de $\gamma_i$ (qui sont disjoints).

\noindent 
3. les $\gamma_i$ sont des géodésiques minimisantes de longueur $l_\h(\gamma_i)=|v_i|_{st}=l_i^{-1}$.

Alors la proposition est vérifiée.
}
\vskip7pt

\noindent {\bf Démonstration du lemme.} On se fixe une métrique lisse $h'$ sur $M$, qui co\"incide avec $\h$ sur les $U_\epsilon(\gamma_i)$. On peut choisir une fonction $\zeta$ lisse, telle que  $\zeta=1$ sur $U_{\epsilon/3}(\gamma_i)$, $\zeta=L$ sur $M \setminus \cup_{i=1}^k U_{2\epsilon/3}(\gamma_i)$ et qui soit croissante sur chaque $U_\epsilon(\gamma_i)$ respectivement à la coordonnée radiale. On choisit $L$ suffisamment grand pour que la métrique $h=\zeta.h'\geq \h$. Pour $i=1,\ldots,k$, $v_i \in \B(h)$. On conclut alors en utilisant le lemme {\bf B1} et la convexité de la boule stable. \hspace{\stretch{1}} $\Box$\\

\medskip
Nous allons maintenant construire une métrique simpliciale satisfaisant les hypothèses du lemme {\bf B4}. Soit $\T^b=H_1(M,\R) / i(H_1(M,\Z))$ le jacobien de $M$ ($i$ désigne l'inclusion canonique). On considère une métrique euclidienne $h_0$ sur $\T^b$ telle que pour toute courbe $\gamma$ telle que $[\gamma]=\sum_{i=1}^k \alpha_i.v_i$ où $\alpha_i \in \N$ pour tout $i=1,\ldots,k$, \linebreak $l_{h_0}(\gamma) \geq \sum_{i=1}^k |\alpha_i|.l_i^{-1}+3$. Fixons un point $x_0$ de $\T^b$. Soient $\vee_{i=1}^{k}C_i$ un bouquet de k cercles pointé en un point $x_1$. On forme un graphe $\Gamma$ en ajoutant à $x_1$ une arête  notée $a$ dont une des extrémités est libre. Soit $h'_0$ une métrique simpliciale sur $\Gamma$ telle que $l_{h'_0}(C_i)=l_i^{-1}$ et $l_{h'_0}(a)=1$. On fixe sur $\T^b$ une triangulation linéaire telle que $x_0$ soit un $0$-simplexe. On considère le complexe simplicial $(P_1,h_1)$ formé en recollant l'extrémité libre de $(\Gamma,h'_0)$ à $(\T^b,h_0)$ en $x_0$ (voir figure 1).

\begin{figure}[h]
\begin{center}
\includegraphics[width=10cm]{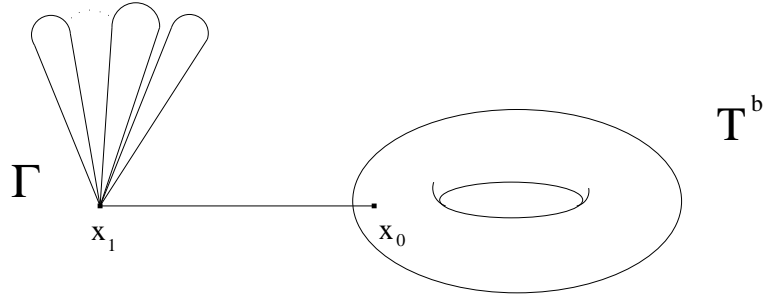}
\caption{Le complexe $P_1$}
\end{center}
\end{figure}

On a une application évidente $f_1 : P_1 \rightarrow \T^b$ qui est l'identité sur $\T^b$ et qui envoie $C_i$ sur $c_i$ la géodésique minimisante réalisant $v_i$ basée en $x_0$ en dilatant les distances.

Comme $f_1$ est l'identité sur $\T_b$, l'application induite en homologie \linebreak
$(f_1)_\ast : H_1(P_1,\R) \rightarrow H_1(\T^b,\R)$ est surjective.  Pour $i=1,\ldots,k$, on réalise le cylindre de l'application $(f_1)_{|C_i} : C_i \rightarrow c_i$ de sorte que l'arête $a$ co\"incide avec une génératrice du cylindre $Z_i$. On obtient ainsi un nouveau complexe simplicial $P_2$ dont $P_1 \subset P_2$ est un sous-complexe. Remarquons également que $T^b \subset P_2$  est un rétract par déformation. On prolonge la métrique $h_1$ de $P_1$ à $P_2$ de la manière suivante. On note $q_i$ la longueur de la concaténation de $\C_i=a \star C_i \star a^{-1}$ et $c_i$.  Soit $p_i : S^1_i \rightarrow \C_i \star c_i $ une application linéaire contractant deux fois les distances. Pour $i=1,\ldots,k$, on considère le cylindre de cette application, et on recolle une demi-sphère deux-dimensionnelle de rayon $q_i / \pi$ isométriquement le long de $S^1_i$ (qui est de longueur $2q_i$). On note $h_2$ la métrique induite par notre construction.

Clairement,
$$
H_1(P_2,\R) \simeq H_1 (\T^b,\R).
$$

Remarquons que pour toute courbe lisse par morceaux $\gamma \subset P_2$, il existe une courbe $\gamma'\subset P_1$ telle que $[\gamma]=[\gamma']$ et $l_{h_2}(\gamma')\leq l_{h_2}(\gamma)$.

\vskip7pt
\noindent
{\bf Lemme B5} {\it $$
\B(h_2)=K.
$$
}
\vskip7pt

\noindent  {\bf Démonstration du lemme.} Il est clair qu'une courbe minimisante de $(P_2,h_2)$ est contenue dans $P_1$. On se fixe une classe entière $v$ appartenant au sous-réseau $\Delta_k$ engendré par $\{v_i\}_{i=1}^k$. Etant donné $n\geq 1$, on choisit $\beta_n$ une courbe réalisant $n.v$. Quitte à modifier notre courbe $\beta_n$, on peut supposer que la courbe est obtenue comme la concaténation de deux courbes : la première $\beta_{n,1}$ est contenue dans $\T^b$ et la seconde $\beta_{n,2}$ dans $\Gamma$. La classe $v$ appartenant au sous-réseau $\Delta_k$, et de la structure de l'homologie réelle unidimensionnelle de $\Gamma$, on déduit que $v$ et $[\beta_{n,2}] \in \Delta_k$ et donc que $[\beta_{n,1}] \in \Delta_k$. On écrit $[\beta_{n,1}]=\sum_{i=1}^k \alpha_i.[C_i]$ et comme $l_{h_0}(\beta_{n,1}) \geq \sum_{i=1}^k |\alpha_i|.l_i^{-1}+3$, on peut diminuer la longueur de $\beta_n$ en rempla\c cant la partie $\beta_{n,1}$ par la courbe $a \star \prod_{i=1}^k C_i^{\alpha_i} \star a^{-1}$. Ceci entraîne que toute courbe minimisante est contenue dans $\Gamma$. Donc $|v|_{st}=\min \{\sum_{i=1}^k |\alpha_i|.l_i^{-1} \mid v=\sum_{i=1}^k \alpha_i.v_i \}$ pour tout $v$ dans le sous-réseau engendré par les $\{v_i\}_{i=1}^k$. Par continuité de la norme et du fait que $rang(\Delta_k)=b$ (comme $int(K) \neq \emptyset$), on en déduit le résultat. \hspace{\stretch{1}} $\Box$\\

\bigskip

On peut remplacer dans notre complexe simplicial $P_2$ les cercles $C_i$ pour $i=1,\ldots,k$ par des tubes $m$-dimensionnels $C_i \times D^{m-1}$ ou $C_i \widehat{\times} D^{m-1}$ (selon l'orientabilité des $U_i$), munis d'une métrique produit (euclidienne sur $D^{m-1}$) dont la projection radiale contracte les distances et telle que les disques $2$-dimensionnels rattachés le soient le long d'un cercle du bord de ce tube. On note à nouveau $(P_2,h_2)$ le complexe simplicial obtenu. Cette opération ne modifie pas la forme de la boule stable.

Soit $\J : M \rightarrow \T^b \subset P_2$. Etant donnés les $k$ voisinages tubulaires topologiques disjoints $U_i$ de $\gamma_i$ pour $i=1,\ldots,k$, on modifie $\J$ dans sa classe d'homotopie de sorte que $\J(\gamma_i)=C_i$ et que $\J(U_i)=C_i \times D^{m-1}$ ou $C_i \widehat{\times} D^{m-1}$.

On peut alors approximer $\J$ par une application $\J_1$ qui soit simpliciale sur $M$ et qui co\"incide avec $\J$ dans chaque $U_i$ pour $i=1,\ldots,k$ (voir \cite{span}). La métrique $\J_1^\ast(h_2)$ tirée en arrière peut être dégénérée en-dehors des voisinages tubulaires. On se fixe une métrique lisse quelconque $h_3$ sur $M$ et $\xi$ une fonction lisse et croissante d'argument $q$ vérifiant
$$\xi(q) = \left\{
\begin{array}{rl}
0 & \mbox{si} \hskip5pt {q \in \cup_{i=1}^k U_i({2\varepsilon \over 3})} ;\\
1 &  \mbox{si} \hskip5pt {q \in M  \setminus \cup_{i=1}^k U_i({\varepsilon})}.
\end{array} \right. $$
 Alors $\J_1 : (M, \J_1^\ast(h_2)+\xi h_3) \rightarrow (P_2,h_2)$ contracte les distances et la métrique ainsi construite vérifie les propriétés annoncées. \hspace{\stretch{1}} $\Box$\\

\bigskip

\subsection{Démonstration du théorème B}

Sur la vari\'et\'e initiale $(M, g)$, on consid\`ere maintenant
les objets construits dans le lemme {\bf B2} et la proposition ci-dessus. 
On choisit deux fonction lisses et croissantes d'argument $t$ vérifiant
$$\phi(t) = \left\{
\begin{array}{rl}
0 & \mbox{si} \hskip5pt t\leq {\varepsilon \over 3}; \\
t &  \mbox{si} \hskip5pt {2\varepsilon \over 3} \leq t;
\end{array} \right. $$
et
$$\psi(t) = \left\{
\begin{array}{rl}
0 & \mbox{si} \hskip5pt t\leq {\varepsilon \over 9} ;\\
1 &  \mbox{si} \hskip5pt {2\varepsilon \over 9} \leq t;
\end{array} \right. $$
o\`u $\varepsilon$ est le rayon des voisinages tubulaires
$\{U_i\}_{i=1}^k$ du lemme {\bf B2}. On consid\`ere un point $q \in U_i$
et ${\bf r}(q)$ la fonction radiale d\'efinie ci-dessus.
Soit $\eta_q(s)$ l'unique g\'eodesique locale passant par $q$, orthogonale
\`a $\gamma_i$ et telle que $\eta_q(0) \in \gamma_i$. 
On d\'efinit les applications
$$ G_i: U_i \longrightarrow U_i ; \hskip5pt 1 \leq i \leq k $$
en posant $G_i(q) =  \eta_q(\phi({\bf r}(q)))$.
Remarquons que $G_i(q) = q$ si ${\bf r}(q) \geq {2\varepsilon \over 3}$ et que $G_i(q)=\eta_q(0)$ ${\bf r}(q) \leq {\varepsilon \over 3}$.  Les applications locales  $\{G_i\}_{i=1}^k$ se prolongent en une
application 
$$ G : M \longrightarrow M .$$
Remarquons que $G$ induit en homologie l'application identit\'e.

On d\'efinit maintenant les fonctions locales suivantes 
$$
\Psi_i : U_i \longrightarrow {\R} ; \hskip5pt 1 \leq i \leq k
$$
en posant $\Psi_i(q) =  \psi({\bf r}(q)))$.
Remarquons que $\Psi_i(q) = 1$ si ${\bf r}(q) \geq {2\varepsilon \over 9}$ et que \linebreak $\Psi_i(q) = 0$ si ${\bf r}(q) \leq {\varepsilon \over 9}$.
Ces fonctions locales  $\{\Psi_i\}_{i=1}^k$ se prolongent en une
fonction $$ \Psi : M \longrightarrow {\R} .$$
Choisissons finalement $B>0$ tel que
$$
 g' = (1 + B \Psi)g_1  \geq g_1 + \Psi G_*(h).
$$
On voit facilement par construction et de la propriété 2. de la proposition {\bf B3} que $G:(M,g') \rightarrow (M,h)$ contracte les distances. 
Du lemme {\bf B1}, on tire
$${\B}(g') \subset {\B}(h) = {K}.$$
Par construction, les m\'etrique $g'$ et $g_1$ co\"{\i}ncident
dans les $\varepsilon \over 9$-voisinages tubulaires des $\{\gamma_i\}_{i=1}^n$.
Ceci avec la propriété 1. du lemme {\bf B2} implique que les sommets de ${K}$
appartiennent \`a ${\B}(g')$. La convexit\'e nous donne alors le r\'esultat.

\end{document}